\documentclass{gen-j-l}
\usepackage{amsmath}
\usepackage{amsfonts}

\newtheorem{theorem}{Theorem}[section]
\newtheorem{lemma}[theorem]{Lemma}
\theoremstyle{definition}
\newtheorem{definition}[theorem]{Definition}

\theoremstyle{remark}
\newtheorem{remark}[theorem]{Remark}
\numberwithin{equation}{section}
\theoremstyle{plain}

\newtheorem{corollary}{Corollary}

\newtheorem{proposition}{Proposition}

\copyrightinfo{2001}{enter name of copyright holder}
\input{tcilatex}

\begin{document}
\title[Hitting time in regular sets and logarithm law for rapidly mixing
systems]{Hitting time in regular sets and logarithm law for rapidly mixing
dynamical systems}
\author{Stefano Galatolo}
\email{s.galatolo@docenti.ing.unipi.it}
\urladdr{http://}
\address{Dipartimento di Matematica Applicata, Universita di Pisa, via
Buonarroti 1 Pisa}

\begin{abstract}
We prove that if a system has superpolynomial (faster than any power law)
decay of correlations (with respect to Lipschitz observables) then the time $%
\tau (x,S_{r})$ needed for a typical point $x$ to enter for the first time a
set $S_{r}=\{x:f(x)\leq r\}$ which is a sublevel of a Lipschitz funcion $f$
scales as $\frac{1}{\mu (S_{r})}$ i.e.%
\begin{equation*}
\underset{r\rightarrow 0}{\lim }\frac{\log \tau (x,S_{r})}{-\log r}=\underset%
{r\rightarrow 0}{\lim }\frac{\log \mu (S_{r})}{\log (r)}.
\end{equation*}%
This generalizes a previous result obtained for balls. We will also consider
relations with the return time distributions, an application to observed
systems and to the geodesic flow of negatively curved manifolds.
\end{abstract}

\thanks{}
\date{23 May 2009}
\maketitle

\section{Introduction and statement of results}

Let $(X,T,\mu )$ be an ergodic system on a metric space $X$ and fix a point $%
x_0 \in X$. For $\mu$-almost every $x \in X$, the orbit of $x$ goes closer
and closer to $x_{0}$ entering (sooner or later) in each positive measure
neighborhood of the target point $x_{0}$.

For several applications it is useful to quantify the speed of approaching
of the orbit of $x$ to $x_{0}$. In the literature this has been done in
several ways, with more or less precise estimations or considering different
kind of target sets.

A\ general approach is to consider a family of sets $S_{r}$ indexed by a
real parameter $r$ containing $x_{0}$ and give an estimation for the time
needed for the orbit of a point $x$ to enter in $S_{r}$%
\begin{equation*}
\tau (x,S_{r})=\min \{n\in \mathbb{N}^{+}:T^{n}(x)\in S_{r}\}.
\end{equation*}

If $X$ is a metric space, the most natural choice is to take $%
S_{r}=B_{r}(x_{0})$ (the ball of radius $r$ ). In this case several
estimations are known for the behavior of $\tau (x,S_{r})$ as $r\rightarrow 0
$. For example if the system has fast decay of correlations or is a circle
rotation, or an interval exchange map with generic arithmetical properties
then for a.e. $x$%
\begin{equation}
\lim_{r\rightarrow 0}\frac{\log \tau (x,B_{r}(x_{0}))}{-\log r}=d_{\mu
}(x_{0})\   \label{1}
\end{equation}%
(see \cite{galatolo,G2,kimseo,KiM}) it is worth to remark that in this case $%
\tau (x,B_{r}(x_{0}))\sim r^{-d_{\mu }}\sim \frac{1}{\mu (S_{r})}$ how it is
natural to expect in a system having "stochastic" behavior: roughly
speaking, as the target shrinks you need more and more iterations to reach
it, but mixing makes geometry of the target to be "forgotten" and only its
measure "recalled", making the time needed to enter, to be inversely
proportional to this measure.

On the other hand it is worth to remark that there are mixing systems
(having particular arithmetical properties and slow decay of correlations)
for which $\lim \inf_{r\rightarrow 0}\frac{\log \tau (x,B_{r}(x_{0}))}{-\log
r}=\infty >d_{\mu }(x_{0})$ (see \cite{GP09} \cite{GRS??}).

In some cases, even if $X$ is a metric space it is interesting to look to
different target sets: for example considering a tubular neighborhood of a
"singular set" and asking how much time we need to approach the singular set
at a distance \thinspace $r$ (see \cite{GD}).

Another interesting situation is when $X$ has a local product structure and
one is interested to approach only some coordinates, for example when $X$ is
a tangent bundle of a manifold $M$ and the dynamics is given by the geodesic
flow, here one can be interested to approach a given target in $M$ giving no
importance to the other coordinates (see \cite{Ma},\cite{Su}, e.g.). Similar
examples are given if one consider the structure given by stable and
unstable directions (\cite{CK}). Another particularly interesting case where
set other than balls become interesting is the case of observed systems (see
section \ref{obs} and \cite{SR} for a similar  point of view in quantitative
recurrence).

We remark that in some of the cited papers the problem of estimating the
speed of approaching a target point is not always directly stated in the
above form, considering the hitting time to a set, but considering the
behavior of some distance to a point or by the so called dynamical Borel
Cantelli lemma (see e.g. \cite{CK},\cite{Dol} which can give a slightly more
precise estimation for typical hitting times in small sets than the one
given in \ref{1}) general relations between these approaches can be found in 
\cite{GP09} and \cite{GK}.

In this note we take the point of view of equation \ref{1} and consider
target sets of the form $S_{r}=\{x\in X~,~f(x)\leq r\}$ where $f$ is a
Lipschitz function. The main result (see theorem \ref{maine} for a precise
statement) is a generalization of the one given in \cite{galatolo} for this
more general family of sets: \emph{if the system has superpolynomial decay
of correlations with respect to Lipschitz observables, then the power law
behavior of the hitting time in regular sets }$S_{r}$\emph{\ as defined in
the abstract satisfies a "logarithm law" of the form}

\begin{equation}
\lim_{r\rightarrow 0}\frac{\log \tau (x,S_{r})}{-\log r}=d(f)\ 
\end{equation}%
where $d(f)=\underset{r\rightarrow 0}{\lim }\frac{\log \mu (S_{r})}{\log (r)}%
,$ generalizing the formula for the local dimension which is in the right
hand of eq. \ref{1}.

We will also give relations limit return time statistic (see Thm. \ref{GAN2}%
) in the sets $S_{r}$ a similar statement for observed systems (see section %
\ref{obs} ).and an application (see section \ref{flowx}) to geodesic flows
on negatively curved manifolds.

\section{Setting and basic results}

Let $f$ be a Borel measurable function such that $f\geq 0$ on $X$. Let
consider sublevel sets $S_{r}=\{x\in X~,~f(x)\leq r\}$. Let us consider the
power law behavior of the hitting time to the set $S_{r}$ as $r\rightarrow 0$

\begin{equation}
\overline{R}(x,f)=\underset{r\rightarrow 0}{\lim \sup }\frac{\log \tau
(x,S_{r})}{-\log (r)},\underline{R}(x,f)=\underset{r\rightarrow 0}{\lim \inf 
}\frac{\log \tau (x,S_{r})}{-\log (r)}.
\end{equation}%
In this way $\tau (x,S_{r})\sim r^{-R(x,f)}$ for small $r$. We remark that
the indicator $R(x,y)$ of \cite{galatolo} is obtained when $f(x)=d(x,y)$ and
the indicator $R(x)$ of \cite{BS} is obtained as a further special case when 
$x=y$.

Let us recall that the definition of local dimension of a measure on a
metric space. If $X$ is a metric space and $\mu $ is a measure on $X$ the
local dimension of $\mu $ at $x$ is defined as $d_{\mu }(x)=%
\mathrel{\mathop{lim}\limits_{r\rightarrow 0}}\frac{log(\mu (B(x,r)))}{log(r)%
}$ (when the limit exists). Conversely, the upper local dimension at $x\in X$
is defined as $\overline{d}_{\mu }(x)=\mathrel{\mathop{limsup}\limits_{r%
\rightarrow 0}}\frac{log(\mu (B(x,r)))}{log(r)}$ and the lower local
dimension $\underline{d}_{\mu }(x)$ is defined in an analogous way by
replacing $limsup$ with $liminf$. If $\overline{d}_{\mu }(x)=\underline{d}%
_{\mu }(x)=d$ almost everywhere the system is called exact dimensional. In
this case many notions of dimension of a measure will coincide. In
particular $d$ is equal to the dimension of the measure: $d=\inf \{\dim
_{H}Z:\mu (Z)=1\}$ (see e.g. \cite{P}). By analogy with the definition of
local dimension let us consider%
\begin{equation}
\overline{d}(f)=\underset{r\rightarrow 0}{\lim \sup }\frac{\log \mu (S_{r})}{%
\log (r)},\underline{d}(f)=\underset{r\rightarrow 0}{\lim \inf }\frac{\log
\mu (S_{r})}{\log (r)}  \label{df}
\end{equation}

Between those two indicators there is a general relation, which follow by a
direct application of the classical Borel Cantelli lemma:

\begin{proposition}
\label{GAN}Let $f$ be as above. Then%
\begin{equation}
\overline{R}(x,f)\geq \overline{d}(f)~,~\underline{R}(x,f)\geq \underline{d}%
(f)
\end{equation}%
$\mu $-a.e.
\end{proposition}

Before to prove the above proposition we point out an elementary remark on
the behavior of real sequences which will be often used in the following.

\begin{lemma}
\label{lemmino} Let $r_{n}$ be a decreasing sequence such that $%
r_{n}\rightarrow 0$. Suppose that there is a constant $c>0$ satisfying $%
r_{n+1}>cr_{n}$ eventually as $n$ increases. Let $\tau_{r}: \mathbb{R}
\rightarrow \mathbb{R}$ be decreasing. Then $\liminf_{n \rightarrow \infty }%
\frac{\log \tau _{r_{n}}}{-\log r_{n}} =\liminf_{r\rightarrow 0}\frac{\log
\tau _{r}}{-\log r}$ and $\limsup_{n\rightarrow \infty }\frac{\log \tau
_{r_{n}}}{-\log r_{n}}=\limsup_{r\rightarrow 0}\frac{\log \tau _{r}}{-\log r}
$.
\end{lemma}

\begin{proof}
(of proposition \ref{GAN}) First we prove\emph{\ $\underline{R}(x,f)\geq 
\underline{d}_{\mu }(f)$. }Let us consider the set of $x$ where $\underline{R%
}(x,f)<\underline{d}_{\mu }(f)$ , we will prove that this set has zero
measure. By the above lemma we will consider a sequence of radii $r_{k}$ of
the form $r_{k}=2^{-k}$. Let us consider $d<\underline{d}_{\mu }(f)$, then
we have eventually that $\mu (S_{2^{-k}})\leq 2^{-dk}$. Let us consider $%
d^{\prime }<d$ and%
\begin{equation*}
A(d^{\prime })=\{x\in X|\underline{R}(x,f)\leq d^{\prime }\}.
\end{equation*}%
By definition of $\underline{R}(x,f)$, it holds that for each $m$%
\begin{equation}
A(d^{\prime })\subset \cup _{n\geq m}\underset{i\leq 2^{n(\frac{d+d^{\prime }%
}{2})}}{\cup }T^{-i}(S_{2^{-n}})
\end{equation}%
but $\mu (\underset{i\leq 2^{n(\frac{d+d^{\prime }}{2})}}{\cup }%
T^{-i}(S_{2^{-n}}))\leq 2^{n(\frac{d+d^{\prime }}{2})}\ast 2^{-dn}$ and $n(%
\frac{d+d^{\prime }}{2})-dn<0$, hence this sequence of sets has summable
measure and by the classical Borel Cantelli lemma $\mu (A(d^{\prime }))=0$,
proving the first inequality.

\emph{Now we prove $\overline{R}(x,f)\geq \overline{d}_{\mu }(f)$.} Again,
we can suppose $r_{k}$ to be of the form $r_{k}=2^{-k}$. Suppose $d^{\prime
}<\overline{d}_{\mu }(f)$, let us consider 
\begin{equation*}
A(d^{\prime })=\{x\in X|\overline{R}(x,f)<d^{\prime }\}.
\end{equation*}%
If $0<d^{\prime }<d<\overline{d}_{\mu }(f)$ then there is a sequence $n_{k}$
such that 
\begin{equation}
\mu (S_{2^{-n_{k}}})<2^{-dn_{k}}\ \mbox{ for each }k.  \label{misu}
\end{equation}%
On the other side for each $x\in A(d^{\prime })$ the relation $\tau
(x,S_{2^{-n}})<2^{\frac{d+d^{\prime }}{2}n}$ must hold eventually. Let us
consider 
\begin{equation}
C(m)=\{x\in A(d^{\prime })|\forall n\geq m,\tau (x,S_{2^{-n}})<2^{\frac{%
d+d^{\prime }}{2}n}\}.
\end{equation}%
This is an increasing sequence of sets \textquotedblleft
converging\textquotedblright\ to $A$. If we prove that $\underset{%
m\rightarrow \infty }{\lim \inf }\mu (C(m))=0$ the statement is proved. By
the definition of $C(m)$ we see that 
\begin{equation}
C(n_{k})\subset \underset{i\leq 2^{\frac{d+d^{\prime }}{2}n_{k}}}{\cup }%
T^{-i}(S_{2^{-n_{k}}})
\end{equation}%
the latter is made of $2^{\frac{d+d^{\prime }}{2}n_{k}}$ sets, whose measure
can be estimated by Eq. \ref{misu}, because $T$ is measure preserving. Then $%
\mu (C(n_{k}))\leq 2^{\frac{d+d^{\prime }}{2}n_{k}}\ast 2^{-dn_{k}}$ and $%
\mu (C(n_{k}))$ goes to $0$ as $k\rightarrow \infty $.
\end{proof}

\section{Fast mixing systems}

As it is well known, in a mixing system we have $\mu (A\cap
T^{-n}(B))\rightarrow \mu (A)\mu (B)$ for each measurable sets $A,B.$ The
speed of convergence of the above limit can be arbitrarily slow (depending
on $T$ but also on the shape of the sets $A,B$). In many systems however the
speed of convergence can be estimated for sets having some regularity.

Let us remark that considering $1_{A}(x)=\left\{ 
\begin{array}{c}
1~if~x\in A \\ 
0~if~x\notin A%
\end{array}%
\right. $ the mixing condition becomes $\int 1_{B}\circ T^{n}1_{A}d\mu
\rightarrow \int 1_{A}d\mu \int 1_{B}d\mu .$

\begin{definition}
\label{sup}Let $\phi ,$ $\psi :X\rightarrow \mathbb{R}$ be Lipschitz
observables on $X$. A system $(X,T,\mu )$ is said to have superpolynomial
decay of correlations with respect to Lipschitz observables, if 
\begin{equation*}
|\int \phi \circ T^{n}\psi d\mu -\int \phi d\mu \int \psi d\mu |\leq
\left\vert \left\vert \phi \right\vert \right\vert \left\vert \left\vert
\psi \right\vert \right\vert \Phi (n)
\end{equation*}%
with $\Phi $ having superpolynolmial decay, i.e. $\lim n^{\alpha }\Phi (n)=0,
$ $\forall \alpha >0.$
\end{definition}

Here $||~||$ is the Lipschitz norm. This is one of the weakest requirements
on the space of observables. Decay of correlations with respect to Holder
observables implies it. In the remaining part of the section we will prove
the following result

\begin{theorem}
\label{maine}If $f:X\rightarrow \mathbb{R}^{+}$ is $\ell $-Lipschitz, the
system has superpolynomial decay of correlations, as above, and $d(f)$
exists, then for a.e. $x$ it holds 
\begin{equation}
R(x,f)=d(f).
\end{equation}
\end{theorem}

Before to prove the theorem we will need some preliminary lemmas: the first
is technical and allow to use decay of correlation to estimate the measure
of certain intersections of sublevels.

\begin{lemma}
\label{uno}Let $r_{n}\rightarrow 0$ and $S_{r_{k}}$ be a sequence of nested
sets as above, let $A_{k}=T^{-k}(S_{r_{k}})$ and let us write $A_{-1}=X$. If 
$(X,T,\mu )$ is a system satisfying definition \ref{sup} then there is $N$
such that when $k>j>N$%
\begin{equation}
\mu (A_{k}\cap A_{j})\leq \mu (A_{k-1})\mu (A_{j-1})+\frac{4\ell ^{2}~\Phi
(k-j)}{(r_{k-1}-r_{k})(r_{j-1}-r_{j})}.  \label{ball}
\end{equation}
\end{lemma}

\begin{proof}
Let 
\begin{equation}
f_{n}(x)=\left\{ 
\begin{array}{c}
\frac{r_{n-1}-f(x)}{r_{n-1}-r_{n}}~if~x\in S_{r_{n-1}}-S_{r_{n}} \\ 
1~if~x\in S_{r_{n}} \\ 
0~if~x\notin S_{r_{n-1}}%
\end{array}%
\right. 
\end{equation}%
these are $\frac{\ell }{r_{n-1}-r_{n}}$-Lipschitz functions, and with
support in $S_{r_{n-1}}.$ Notice that $\mu (S_{r_{n}})\leq \int f_{n}(x)d\mu
\leq \mu (S_{r_{n-1}})$. The Lipschitz norm is such that $%
||f_{k}||_{Lip}\leq \frac{\ell }{r_{n-1}-r_{n}}+1$. If $N$ is big enough
that $r_{n-1}-r_{n}\leq \ell $, $\forall n\geq N,$ then $||f_{k}||_{Lip}\leq 
\frac{2\ell }{r_{n-1}-r_{n}}$. Let $k>j>N$. Since $\mu $ is preserved%
\begin{equation*}
\mu (A_{k}\cap A_{j})=\mu (T^{-k+j}(S_{r_{k}})\cap S_{r_{j}})\leq \int
f_{k}\circ T^{k-j}~f_{j}~d\mu .
\end{equation*}%
By decay of correlations%
\begin{eqnarray*}
\int f_{k}\circ T^{k-j}~f_{j}~d\mu  &\leq &\int f_{k}d\mu \int f_{j}d\mu
+\left\vert \left\vert f_{k}\right\vert \right\vert _{Lip}\left\vert
\left\vert f_{j}\right\vert \right\vert _{Lip}\Phi (k-j)\leq  \\
&\leq &\mu (A_{k-1})\mu (A_{j-1})+\left\vert \left\vert f_{k}\right\vert
\right\vert _{Lip}\left\vert \left\vert f_{j}\right\vert \right\vert
_{Lip}\Phi (k-j)
\end{eqnarray*}%
which gives the statement.
\end{proof}

The second Lemma we will use is a sort of dynamical Borel Cantelli lemma for
systems having fast decay of correlations ( \cite{galatolo} , Lemma 7)

\begin{lemma}
\label{due}Let $S_{k}$ be a decreasing sequence of measurable sets such that 
\begin{equation*}
\underset{k\rightarrow \infty }{\lim \inf }\frac{\log (\sum_{0}^{k}\mu
(S_{k}))}{\log (k)}=z>0.
\end{equation*}%
Let $A_{k}=T^{-k}(S_{k})$ and let us suppose that the system is such that
when $k>j$%
\begin{equation}
\mu (A_{k}\cap A_{j})\leq \mu (A_{k-1})\mu (A_{j-1})+k^{c_{1}}j^{c_{2}}\Phi
(k-j)  \label{mixxx}
\end{equation}%
with $\Phi $ having superpolynomial decay and $c_{1},c_{2}\geq 0$. Then
posing 
\begin{equation}
Z_{k}(x)=\sum_{0}^{k}1_{A_{i}}(x)
\end{equation}
we have $\frac{Z_{k}}{E(Z_{k})}\rightarrow 1$ in the $L^{2}$ norm and almost
everywhere.
\end{lemma}

We remark that, since the set sequence is decreasing the condition $\mu
(A_{k}\cap A_{j})\leq \mu (A_{k-1})\mu (A_{j-1})+...$ is slightly more
relaxed than $\mu (A_{k}\cap A_{j})\leq \mu (A_{k})\mu (A_{j})+...$. This is
a technical point which will allow us to use Lemma \ref{uno} without further
technical complications\ and apply the Lemma \ref{due} \ to the sequence $%
S_{r_{k}}$. This will allow to prove the main result:

\begin{proof}
\emph{(of Thm. \ref{maine})} Let us prove $\overline{R}(x,f)\leq d_{\mu }(f)$
for almost each $x.$ For simplicity of notations let us set $d=d_{\mu }(f)$.
We recall that this implies \underline{$R$}$(x,f)\leq d$ and the opposite
inequalities come from theorem \ref{GAN}. Let us consider $0<\beta <\frac{1}{%
d_{\mu }(f)}$, the sequence $r_{k}=k^{-\beta }$ and set $S_{k}=S_{r_{k}}$
(we remark that if the result is proved for such a subsequence, then it
holds for all subsequences, see lemma \ref{lemmino}). For each small $%
\epsilon <\beta ^{-1}-d$, eventually $\mu (S_{k})\geq (r_{k})^{d+\epsilon
}=k^{-\beta (d+\epsilon )}$ and if $k$ is big enough $\sum_{0}^{k}\mu
(S_{k})\geq Ck^{1-\beta (d+\epsilon )}.$ Since $\epsilon $ is arbitrary we
have 
\begin{equation*}
\underset{k\rightarrow \infty }{\lim \inf }\frac{\log (\sum_{0}^{k}\mu
(S_{k}))}{\log (k)}\geq 1-\beta (d+\epsilon )>0.
\end{equation*}%
Moreover, $r_{k-1}-r_{k}\sim k^{-\beta -1}.$ Hence we can apply Lemma \ref%
{uno} and \ref{due} to the sequence $S_{k}$ and obtain that for such a
sequence $\underset{n\rightarrow \infty }{lim}\frac{Z_{n}}{E(Z_{n})}=1,$ $%
\mu -$almost everywhere.

Let us now consider $\epsilon ^{\prime }>0$ such that $\beta (d+\epsilon
^{\prime })>1$ for $\beta $ as above, near to $\frac{1}{d}$. Moreover let us
consider $\varepsilon >0$ so small that $\beta (d+\varepsilon )<1$ and $%
\beta (d+\epsilon ^{\prime })-\frac{1-\beta (d-\varepsilon )}{1-\beta
(d+\varepsilon )}>0$. Let us consider $x$ such that $\overline{R}%
(x,f)>d+\epsilon ^{\prime }$, then for infinitely many $n$, $\tau
(x,S_{n})>n^{\beta (d+\epsilon ^{\prime })}$. Then 
\begin{equation*}
x\notin \cup _{0\leq i\leq \left\lfloor n^{\beta (d+\epsilon ^{\prime
})}\right\rfloor }T^{-i}(S_{_{n}})
\end{equation*}%
and in particular 
\begin{equation*}
x\notin \cup _{n\leq i\leq \left\lfloor n^{\beta (d+\epsilon ^{\prime
})}\right\rfloor }T^{-i}(S_{_{n}})\supset \cup _{n\leq i\leq \left\lfloor
n^{\beta (d+\epsilon ^{\prime })}\right\rfloor }T^{-i}(S_{_{i}}),
\end{equation*}%
which implies that there is a sequence $n_{i}$ such that $%
Z_{n_{i}}(x)=Z_{\left\lfloor n_{i}^{\beta (d+\epsilon ^{\prime
})}\right\rfloor }(x)$ ($Z$ was defined in lemma \ref{due}) for each $i.$
Now let us consider $E(Z_{n_{i}})$ and $E(Z_{\left\lfloor n_{i}^{\beta
(d+\epsilon ^{\prime })}\right\rfloor }).$ By the definition of $d=d_{\mu
}(f)$,\ when $i$ is big enough 
\begin{equation*}
i^{-\beta (d+\varepsilon )}<\mu (S_{i})<i^{-\beta (d-\varepsilon )}
\end{equation*}%
then there are constants $k_{1}$ and $k_{2}$ such that when $n$ is big
enough $k_{1}n^{1-\beta (d+\varepsilon )}<E(Z_{n})<k_{2}n^{1-\beta
(d-\varepsilon )}$. From this we have that if $i$ is big enough%
\begin{eqnarray*}
\frac{E(Z_{n_{i}})}{E(Z_{\left\lfloor n_{i}^{\beta (d+\epsilon ^{\prime
})}\right\rfloor })} &\leq &\frac{k_{2}n_{i}^{1-\beta (d-\varepsilon )}}{%
k_{1}\left\lfloor n_{i}^{\beta (d+\epsilon ^{\prime })}\right\rfloor
^{(1-\beta (d+\varepsilon ))}}\sim  \\
&\sim &\frac{k_{2}}{k_{1}}n_{i}^{(1-\beta (d-\varepsilon ))-\beta
(d+\epsilon ^{\prime })(1-\beta (d+\varepsilon ))}.
\end{eqnarray*}%
By the assumptions on $\varepsilon ,$ $(1-\beta (d-\varepsilon ))-\beta
(d+\epsilon ^{\prime })(1-\beta (d+\varepsilon ))=(1-\beta (d+\varepsilon ))(%
\frac{1-\beta (d-\varepsilon )}{1-\beta (d+\varepsilon )}-\beta (d+\epsilon
^{\prime }))<0,$ hence%
\begin{equation*}
\underset{i\rightarrow \infty }{\lim }\frac{E(Z_{n_{i}})}{E(Z_{\left\lfloor
n_{i}^{\beta (d+\epsilon ^{\prime })}\right\rfloor })}=0.
\end{equation*}%
Since $n_{i}$ was chosen such that $Z_{n_{i}}(x)=Z_{\left\lfloor
n_{i}^{\beta (\overline{d}_{\mu }(f)+\epsilon ^{\prime })}\right\rfloor }(x)$
this implies that 
\begin{equation}
\frac{Z_{n_{i}}(x)}{E(Z_{n_{i}})}\frac{E(Z_{\left\lfloor n_{i}^{\beta
(d+\epsilon ^{\prime })}\right\rfloor })}{Z_{\left\lfloor n_{i}^{\beta
(d+\epsilon ^{\prime })}\right\rfloor }(x)}=\frac{E(Z_{\left\lfloor
n_{i}^{\beta (d+\epsilon ^{\prime })}\right\rfloor })}{E(Z_{n_{i}})}%
\rightarrow \infty 
\end{equation}
as $i$ increases. Then is not possible that $\underset{n\rightarrow \infty }{%
\lim }\frac{Z_{n}(x)}{E(Z_{n}(x))}=1.$ This, implies that $\overline{R}%
(x,f)>d+\epsilon ^{\prime }$ on a zero measure set. Finally, since $\epsilon
^{\prime }$ can be chosen to be arbitrarily small we have the statement.
\end{proof}

\subsubsection{Geodesic flow and its time one map\label{flowx}}

As a simple application of the main theorem (\ref{maine}) we give a sort of
logarithm law for the geodesic flow in a negatively curved manifold, which
is a sort of discrete time version of the main result of \cite{Ma}. We
remark that since we only need an estimation for the decay of correlation,
we can apply the deep results available for this kind of flows and consider
manifolds with (variable) negative curvature instead of constant curvature
ones. We will use the following result from \cite{Li}:

\begin{theorem}
The geodetic flow $T^{t}$ of a $C^{4}$ manifold with strictly negative
curvature is exponentially mixing with respect to Holder observables: there
exists $C$ and $\sigma >0$ such that 
\begin{equation}
|\int \phi \circ T^{t}\psi d\mu -\int \phi d\mu \int \psi d\mu |\leq
C\left\vert \left\vert \phi \right\vert \right\vert _{\alpha }\left\vert
\left\vert \psi \right\vert \right\vert _{\alpha }e^{-\sigma t}.
\end{equation}
\end{theorem}

Theorem \ref{maine} hence gives

\begin{proposition}
Let $M$ be a $C^{4}$ manifold of dimension $d$ with strictly negative
curvature and $T^{1}M$ be its unitary tangent bundle. Let $\pi
:T^{1}M\rightarrow M$ the canonical projection. If $T$ is the time $1$ map
of the geodesic flow,$\ \mu $ the Liouville measure on $T^{1}M$ , and $d$
the Riemannian distance on $M,$ then for each $p\in M$:%
\begin{equation}
\underset{n\rightarrow \infty }{\lim \sup }\frac{-\log d(p,\pi (T^{n}x))}{%
\log n}=\frac{1}{d}
\end{equation}%
holds for almost each $x\in T^{1}M.$
\end{proposition}

\begin{proof}
By the above theorem, the time $1$ map associated to the flow has
exponential decay of correlations for Holder observables and hence for
Lipschitz too. Let us consider the function $f:T^{1}M\rightarrow \mathbb{R}$
given by $f(x)=d(\pi (x),p)$. Since the Liouville measure is absolutely
continuous, with density bounded away from zero, the associated sets $S_{r}$
are such that $\frac{\log \mu (S_{r})}{\log (r)}\rightarrow d$, then by
Theorem \ref{maine}%
\begin{equation}
\underset{r\rightarrow 0}{\lim }\frac{\log \tau (x,S_{r})}{-\log (r)}=d.
\end{equation}%
Let $d_{n}(x,p)=\min_{i\leq n}\mathrm{dist}(\pi (T^{i}(x)),p)$, let us
suppose that for $n$ big enough $d_{n}=n^{-\alpha (n)}$. Since $\underset{%
n\rightarrow 0}{\lim }\frac{\log \tau (x,S_{n^{-\alpha (n)}})}{-\log
(n^{-\alpha (n)})}=d$ then $\forall \epsilon $ if $n$ big enough%
\begin{equation}
n^{\alpha (n)d-\alpha (n)\epsilon }\leq \tau (x,S_{n^{-\alpha (n)}})
\end{equation}%
but $\tau (x,S_{n^{-\alpha (n)}})\leq n$ hence $n^{\alpha (n)d-\alpha
(n)\epsilon }\leq n$ , $\alpha (n)d-\alpha (n)\epsilon \leq 1$ and then $%
\alpha (n)\leq \frac{1}{d-\epsilon }$. This implies that $\underset{%
n\rightarrow \infty }{\lim \sup }\frac{-\log d_{n}(p,\pi (T^{n}x))}{\log n}%
\leq \frac{1}{d}$. On the other hand there are infinitely many $n$ such that 
$\tau (x,S_{n^{-\alpha (n)}})=n$. This with the same calculation as above
implies the statement.
\end{proof}

\subsection{Observed systems\label{obs}}

An important case where hitting time for sets different than balls become
interesting and our approach very natural is the case of observed systems.
Here the behavior of the system $(X,d,T,\mu )$ is observed trough a
measurable function $F:X\rightarrow Y$ and we look to the time needed for $%
F(T^{n}x)$ to approach $F(x_{0})$ (in \cite{SR} something similar was done
for quantitative recurrence indicators).

The function naturally induces a measure $F^{\ast }(\mu )$ on $Y$, defined
as $F^{\ast }(\mu )[A]=\mu (F^{-1}(A))$ for each measurable set $A\subseteq
Y $. We can then consider the local dimension of the induced measure $%
d_{F^{\ast }(\mu )}$ of $F^{\ast }(\mu )$ \ and the observed hitting times:

\begin{equation}
\tau _{r}^{F}(x,x_{0})=\min \{k\in \mathbb{N}^{+},F(T^{k}(x))\in
B_{r}(F(x_{0}))\}.
\end{equation}

\begin{proposition}
Let us consider $f:X\rightarrow \mathbb{R}$ defined by $%
f(x)=d(F(x),F(x_{0})) $, then for each $x\in X$%
\begin{eqnarray}
d_{F^{\ast }(\mu )}(F(x_{0})) &=&d(f), \\
\tau _{r}^{F}(x,x_{0}) &=&\tau (x,S_{r})
\end{eqnarray}%
where $S_{r}=\{x\in X~,~f(x)\leq r\}$ as in \ref{df} and $d(f)$ is also
defined as in \ref{df}.
\end{proposition}

\begin{proof}
The proof is straightforward from definitions

\begin{eqnarray*}
d_{F^{\ast }(\mu )}(F(x_{0})) &=&\lim_{r\rightarrow 0}\frac{\log (F^{\ast
}(\mu )[B_{r}(F(x_{0}))])}{\log r}=\lim_{r\rightarrow 0}\frac{\log (\mu
(F^{-1}(B_{r}(F(x_{0}))))}{\log r}= \\
&=&\underset{r\rightarrow 0}{\lim }\frac{\log \mu (\{x\in
X~,~d(F(x),F(x_{0}))\leq r\})}{\log (r)}=\underset{r\rightarrow 0}{\lim }%
\frac{\log \mu (S_{r})}{\log (r)}=d(f)
\end{eqnarray*}%
moreover%
\begin{equation*}
\tau _{r}^{F}(x,x_{0})=\min \{k\in \mathbb{N},F(T^{k}(x))\in
B_{r}(F(x_{0}))\}=\tau (x,S_{r}).
\end{equation*}
\end{proof}

This gives the following corollary of the main theorem \ref{maine} for
observed systems

\begin{corollary}
If we consider an observed system as above and: $F:X\rightarrow Y$ is
Lipschitz, the system has superpolynomial decay of correlatinos, as above, $%
d_{F^{\ast }(\mu )}(F(x_{0}))$ exists, then 
\begin{equation*}
\lim_{r\rightarrow 0}\frac{\log \tau _{r}^{F}(x,x_{0})}{-\log (r)}%
=d_{F^{\ast }(\mu )}(F(x_{0}))
\end{equation*}%
$\mu $-a.e.
\end{corollary}

The following theorem form \cite{SR} ensures the existence of the local
dimension for a class of interesting examples of observed systems.

\begin{theorem}
Let $F:\mathbb{R}^{m}\mathbb{\rightarrow }\mathbb{R}^{n}$ be a $C^{\infty }$%
function, let $\mu $ be an absolutely continuous measure on $\mathbb{R}^{m}$%
, then $d_{F^{\ast }(\mu )}$ exists and belongs to the set $\{0,1,...,\min
(m,n)\}$ almost everywhere. More precisely, $d_{F^{\ast }(\mu
)}(f(x))=rank(d_{x}F)$ for $\mu -$almost every $x\in \mathbb{R}^{M}.$
\end{theorem}

\section{A relation with return time statistics}

The return time statistics is a widely studied feature of dynamics (see e.g. 
\cite{LHV}, \cite{BS} and references therein) and has links with other
subjects as the extreme value theory (\cite{F}) Let us consider a measurable
function $f$, the above sets $S_{r}=\{x:f(x)\leq r\}$ and the statistical
distribution of return times in these sets. We say that the return time
statistics of $(X,T)$ converges to $g$ for the sets $S_{r}$, if%
\begin{equation}
\underset{r\rightarrow 0}{\lim }\frac{\mu (\{x\in S_{r},\tau (x,S_{r})\geq 
\frac{t}{\mu (S_{r})}\})}{\mu (S_{r})}=g(t).  \label{1221}
\end{equation}

If $g(t)=e^{-t}$ we say that the system has an exponential return time limit
statistics. Such statistics can be found in several systems with some
hyperbolic behavior in some class of decreasing sets, however other limit
distributions are possible. The following shows that if the logarithm law
does not hold then the return time statistic has a trivial limit.

\begin{theorem}
If $(X,T,\mu )$ is ergodic, the measure is finite, nonatomic and 
\begin{equation}
\underline{R}(x,f)>\overline{d}_{\mu }(f)
\end{equation}%
a.e., then the system has trivial return time statistic in the sets $S_{r}$.
That is, the limit in \ref{1221} exists for each $t$ and $g(t)=0$.
\end{theorem}

Let us consider the set%
\begin{equation}
C_{l,r}=\{x\in S_{r},\tau (x,S_{r})>l\mu (S_{r})^{-1}\}.
\end{equation}%
We remark that if the return time statistic converges to $g$ as above then $%
\underset{r\rightarrow 0}{\lim }\frac{\mu (C_{l,r})}{\mu (S_{r})}=g(l)$. The
above theorem is implied by the following

\begin{lemma}
If there is an $l>0$ such that $\underset{r\rightarrow 0}{\lim \sup }\frac{%
\mu (C_{l,r})}{\mu (S_{r})}>0$ then $\underline{R}(x,f)\leq \overline{d}%
_{\mu }(f).$
\end{lemma}

\begin{proof}
Let us consider a number $l$, a sequence $r_{n}\rightarrow 0$ such that $%
\lim_{n\rightarrow \infty }\frac{\mu (C_{l,r_{n}})}{\mu (S_{r_{n}})}>0$ and
the sets $T^{-1}(C_{l,r_{n}}),T^{-2}(C_{l,r_{n}}),...,T^{-\left\lfloor l\
(\mu (S_{r_{n}})^{-1})\right\rfloor }(C_{l,r_{n}})$. All these sets are
disjoint because if there was $x\in T^{-i}(C_{l,r_{n}})\cap
T^{-j}(C_{l,r_{n}})$ with $i,j\leq \left\lfloor l\ (\mu
(S_{r_{n}})^{-1})\right\rfloor $ then $\tau (T^{min(i,j)}(x),S_{r_{n}})\leq
|i-j|$ and by definition of $C_{l,r_{n}}$, $T^{min(i,j)}(x)$ cannot be
contained in $C_{l,r_{n}}$, leading to a contradiction. The set $%
U_{n}=T^{-1}(C_{l,r_{n}})\cup T^{-2}(C_{l,r_{n}})\cup ...\cup
T^{-\left\lfloor l\ (\mu (S_{r_{n}})^{-1})\right\rfloor }(C_{l,r_{n}})$ is
then such that $\mu (U_{n})\geq \left\lfloor l\ (\mu
(S_{r_{n}})^{-1})\right\rfloor \mu (C_{l,r_{n}})$. Since $\lim_{n}\frac{\mu
(C_{l,r_{n}})}{\mu (S_{r_{n}})}>0$ there is a $c>0$ s.t. $\mu
(C_{l,r_{n}})\geq c\mu (S_{r_{n}})$ for $n$ big enough, then $\mu
(U_{n})\geq c\frac{\left\lfloor l\ (\mu (S_{r_{n}})^{-1})\right\rfloor }{\mu
(S_{r_{n}})^{-1}}>c^{\prime }>0$. We remark that the set of points%
\begin{equation*}
G=\{x\in X\ s.t.\ x\in U_{r_{n}}\ for\ infinitely\ many\ n\}
\end{equation*}%
has positive measure. We remark that if a point $x$ is contained in $%
U_{r_{n}}$ for some $n$ then $T^{i}(x)\in S_{r_{n}}$ for $i\leq \left\lfloor
l\ (\mu (S_{r_{n}})^{-1})\right\rfloor $ and then $\tau (x,S_{r_{n}})\leq $ $%
\left\lfloor l\ (\mu (S_{r_{n}})^{-1})\right\rfloor .$

If $\overline{d}_{\mu }(f)=d$ then for each $\delta >0$, if $n$ is big
enough, $\mu (S_{r_{n}})\geq {r_{n}}^{d+\delta }$. Hence $\tau
(x,S_{r_{n}})\leq $ $\left\lfloor l\ (\mu (S_{r_{n}})^{-1})\right\rfloor
\leq \left\lfloor l~{r_{n}}^{-(d+\delta )}\right\rfloor $ , giving $\frac{%
-log(\tau (x,S_{r_{n}}))}{log(r_{n})}\leq d+\delta $ for infinitely many $n$
then $\underline{R}(x,f)\leq d$. This is true for $x$ belonging in a
positive measure set $G$. Since $\mu $ is ergodic and $R(x,f)=R(T(x),f)$, we
have the statement.
\end{proof}

\begin{remark}
The mixing system without logarithm law given in \cite{GP09} hence has
trivial return limit statistic in each centered sequence of balls, this give
an example of a smooth mixing systems with trivial return time statistics.
\end{remark}


\begin{thebibliography}{99}
\bibitem{BS} Barreira L., Saussol B.,\ \emph{Hausdorff dimension of measures
via Poincar\'{e} recurrence}, Commun. Math. Phys. \textbf{219} (2001),
443--463.

\bibitem{CK} N Chernov, D Kleinbock \emph{Dynamical Borel-Cantelli lemmas
for Gibbs measures }Israel Journal of Mathematics,Vol 122, N 1 ,1-27 (2001)

\bibitem{Dol} Dolgopyat D., \emph{Limit theorems for partially hyperbolic
systems}, Trans. Amer. Math. Soc. \textbf{356} (2004), 1637--1689.

\bibitem{GD} Degli Esposti M.Galatolo S. \emph{Recurrence near given sets
and the complexity of the Casati--Prosen map }Chaos, Solitons and Fractals
Vol. 23,4, 1275-1284 (2005)

\bibitem{F} A.C.M. Freitas, J.M. Freitas, M. Todd, \emph{Hitting Time
Statistics and Extreme Value Theory}, Preprint arXiv:0804.2887, 2008.

\bibitem{galatolo} Galatolo S., \emph{Dimension and waiting time in rapidly
mixing systems}, Math Res. Lett. (2007).

\bibitem{GP09} Galatolo S., Peterlongo P. \emph{Long hitting time, slow
decay of correlations and arithmetical properties.}Preprint arXiv:0801.3109

\bibitem{GRS??} Galatolo S, Rousseau J, Saussol B \emph{Work in preparation }

\bibitem{G} Galatolo S., \emph{Dimension via waiting time and recurrence},
Math. Res. Lett. \textbf{12} (2005), 377--386.

\bibitem{G2} Galatolo S., \emph{\ Hitting time and dimension in axiom A
systems, generic interval exchanges and an application to Birkoff sums.} J.
Stat. Phys. \textbf{123} (2006), 111--124.

\bibitem{GK} Galatolo S., Kim D. H., \emph{The dynamical Borel-Cantelli
lemma and the waiting time problems}, Preprint Arxiv: math.DS/0610213.

\bibitem{HV} Hill R., Velani S., \emph{The ergodic theory of shrinking
targets} Inv. Math. \textbf{119} (1995), 175--198.

\bibitem{kimseo} Kim D. H. and Seo B. K., \emph{The waiting time for
irrational rotations}, Nonlinearity \textbf{16} (2003), 1861--1868.

\bibitem{KM} Kleinbock D. Y. , Margulis G. A., \emph{Logarithm laws for
flows on homogeneous spaces. } Inv. Math. \textbf{138} (1999), 451--494.
265--326

\bibitem{KiM} Kim D H, Marmi S \emph{The recurrence time for interval
exchange maps} Nonlinearity 21 2201-2210 (2008)

\bibitem{LHV} Lacroix Y., Haydn N., Vaienti S., \emph{Hitting and return
times in ergodic dynamical systems}, Ann. Probab. 33 (2005), no. 5,
2043--2050.

\bibitem{Li} Liverani C \emph{On contact Anosov flows, }Ann. Math.(2004)

\bibitem{Ma} Maucourant F. \emph{Dynamical Borel-Cantelli Lemma for
hyperbolic spaces} Israel Journal of mathematics, 152 (2006), p 143-155

\bibitem{P} Pesin Y \emph{Dimension theory in dynamical systems,} Chicago
lectures in Mathematics (1997).

\bibitem{SR} J Rousseau, B Saussol \emph{Poincar\'{e} recurrence for
observations }Arxiv preprint arXiv:0807.0970,

\bibitem{Su} Sullivan D., \emph{Disjoint spheres, approximation by imaginary
quadratic numbers, and the logarithm law for geodesics} Acta Mathematica 
\textbf{149} (1982), 215--237
\end{thebibliography}
\end{document}